\documentclass[twoside,11pt]{article}

\usepackage{amsmath,amsfonts}
\usepackage{ctex}

\newcommand{\D}{\displaystyle}

\newcommand{\dis}{\displaystyle}

\def \qed {\hfill \vrule height6pt width 6pt depth 0pt}

\begin{document}
\title{\bf {Decay of a model system of radiating gas}}

\author{W{\sc enjun} W{\sc ang}\\
\small {\it College of Science, University of Shanghai for Science
and Technology,}\\
\small {\it Shanghai 200093, P.R. China}\\
 \and  Z{\sc higang}
W{\sc u}\thanks{Corresponding author. E-mail: mathzgwu@yahoo.com.cn}\\
\small{\it Department of Mathematics, Hangzhou Normal University,}\\
\small{\it Hangzhou 310036, P.R. China} }

  \date{}

\maketitle

\textbf{{\bf Abstract:}} This paper is concerned with optimal
time-decay estimates of solutions of the Cauchy problem to a model
system of the radiating gas in $\mathbb{R}^n$. Compared to Liu and
Kawashima (2011) \cite{Liu1} and Wang and Wang (2009) \cite{Wang},
without smallness assumption of initial perturbation in $L^1$-norm,
we study large time behavior of small amplitude classical solutions
to the Cauchy problem. The optimal $H^N$-norm time-decay rates of
the solutions in $\mathbb{R}^n$ with $1\leq n\leq4$ are obtained by
applying the Fourier splitting method introduced in Schonbek (1980)
\cite{Schonbek1} with a slight modification and an energy method.
Furthermore, basing on a refined pure energy method introduced in
Guo and Wang \cite{Guo} (2011), we give optimal
$L^p$-$L^2(\mathbb{R}^3)$ decay estimates of the derivatives of
solutions when initial perturbation is bounded in $L^p$-norm with
some $p\in(1,2]$.

\bigbreak \textbf{{\bf Key Words}:} the radiating gas; Fourier
splitting method; negative Sobolev space; optimal time-decay rate.

\section*{ 1.\ Introduction}

In the present paper, we consider the large time behavior of the
Cauchy problem for a model system of radiating gas, taking the form
of
$$
\left\{\begin{array}{l}  u_{t}+\sum\limits_{j=1}^nf_j(u)_{x_j}+{\rm div}q=0,\ \ \ \  x\in \mathbb{R}^n,\ t>0,\\[3mm]
-\nabla{\rm div}q+q+\nabla u=0,\ \ \ \ \ \ \ \ \ \ x\in
\mathbb{R}^n,\ t>0,
 \end{array}
        \right.
        \eqno(1.1)
$$
with initial data
$$
u(x,0)=u_0(x),\ \ \ \  x\in \mathbb{R}^n. \eqno(1.2)
$$
Here unknown functions
$u=u(x,t):\mathbb{R}^n\times[0,\infty)\rightarrow \mathbb{R}$ and
$q=q(x,t):\mathbb{R}^n\times[0,\infty)\rightarrow\mathbb{R}^n$
represent the velocity and radiating heat flux of the gas,
respectively. The notations $\nabla$ and ${\rm div}$ are the
$n$-dimensional gradient and divergence.
$f(u)=(f_1(u),\cdots,f_n(u))\in \mathbb{R}^n$ is a given smooth
function of $u$ satisfying $f_j(u)=O(|u|^2)\ (j=1,\cdots,n)$ for
$u\rightarrow0$.

System (1.1) simplifies the model for the motion of radiating gas in
several space variables. Indeed, in a certain physical situation,
(1.1) is well approximated to the fundamental system describing the
motion of radiating gas:
$$
\left\{\begin{array}{l}  \rho_{t}+{\rm div}(\rho u)=0,\\[3mm]
(\rho u)_t+{\rm div}(\rho u\otimes u+pI)=0\\[3mm]
[\rho(e+\frac{|u|^2}{2})]_t+{\rm div}[\rho
u(e+\frac{|u|^2}{2})+pu+q]=0\\[3mm]
 -\nabla{\rm div}q+a_1q+a_2\nabla\theta^4=0,
 \end{array}
        \right.
        \eqno(1.3)
$$
where $\rho,u,p,e,\theta$ are respectively the mass density,
velocity, pressure, internal energy and absolute temperature of the
gas, while $q$ is the radiative heat flux. $a_1$ and $a_2$ are given
positive constants depending on the gas itself. The simplified model
(1.1) was first investigated by Hamer \cite{Hamer}, and the
reduction of the full system (1.3) to (1.1) was given in \cite{Gao2,
Hamer, Kawashima7}.

A lot of important works have been done on system (1.1). For
one-dimensional case, we refer to \cite{Francesco2,Serre1,Serre2}
for $L^1$ stability results, \cite{Kawashima6,Lattanzio1,Laurenot1}
for a singular limit and relaxation limit,
\cite{Kawashima3,Kawashima4,Kawashima5,Lattanzio2,
Lattanzio3,Lin1,Lin2,LiuH,Nishibata,Nguyen} for shock waves,
\cite{Iguchi,Kawashima1,Kawashima2} for diffusion waves and
\cite{Kawashima7} for rarefaction waves.

However, there are fewer studies for system (1.1) in the case of
multi-dimensional space. Recently, Francesco in \cite{Francesco1}
obtained the global existence of weak entropy solutions of system
(1.1) and the relaxation limits. Later, Wang and Wang in \cite{Wang}
investigated pointwise estimates of solutions to (1.1)-(1.2) by
means of a detailed analysis for Green's function. More recently,
the asymptotic decay rates toward the planar rarefaction waves based
on the $L^2$-energy method are obtained in \cite{Gao1} for
2-dimensions and in \cite{Gao2} for $n$-dimensions ($n=3,4,5)$,
respectively. The asymptotic behavior of solutions to the diffusion
waves was studied in \cite{Liu1,Ruan}. For the related study of
decay rates to the problem (1.1)-(1.2), we mention \cite{Liu1} and
\cite{Duan1} only. \cite{Liu1} studied large time behavior of
solutions to the problem (1.1)-(1.2) with small initial data in
$L^1$-norm perturbation by using a time-weighted energy method. In
\cite{Duan1}, under the assumption that $\|u_0\|_{L^1}$ is bounded,
Duan $et\ al.$ showed that the optimal $L^2$-norm time-decay rate of
solutions is $(1+t)^{-\frac{n}{4}}$, whereas time-decay estimates of
the derivatives of the solutions have not been considered.

The purpose of this paper is to establish the optimal time decay
rates of the solutions and the derivatives of the solutions to the
Cauchy problem (1.1)-(1.2) without smallness assumptions on the
initial data in $L^1$-norm. It seems that the usual energy method
basing on the linearization analysis doesn't work. By the way, our
first decay result is inspired from Schonbek
\cite{Schonbek2,Schonbek3}, where the well-known Fourier splitting
method is established to get optimal decay rate of solutions for the
incompressible Navier-Stokes
 equations in $L^2$-norm or $H^s$-norm. In the present paper,
 we generalize the Fourier splitting method with a slight
modification to deal with the problem (1.1)-(1.2) for $\mathbb{R}^n$
with $1\leq n\leq 4$. Our second decay result is illuminated by a
recent work of Guo and Wang \cite{Guo}, where they developed a new
method to establish optimal time decay rates of solutions to the
Cauchy problem for the compressible Navier-Stokes equations and the
Boltzmann equation. The main idea was to combine scaled energy
estimates with an interpolation between negative and positive
Sobolev norms to get the time decay rate for dissipative equations.
By employing this new method, we obtain the optimal
$L^p$-$L^2(\mathbb{R}^3)$ decay estimates of solutions to the
problem (1.1)-(1.2).

\bigbreak

First, we cite a temporal global existence result established in
\cite{Liu1,Wang}.

\vspace{2mm}

 \noindent\textbf{Proposition 1.1.}(\cite{Liu1}, Theorem 2.1;
\cite{Wang}, Theorem 1.1) Assume that $u_0(x)\in H^N(\mathbb{R}^n)\
(n\geq1)$ for an integer $N\geq[\frac{n}{2}]+2$. There exists a
small positive constant $\delta_0$ such that if
$E_0=\|u_0\|_{H^N(\mathbb{R}^n)}\leq\delta_0$, then the problem
(1.1)-(1.2) has a unique global solution $(u,q)(x,t)$ satisfying
$$\arraycolsep=1.5pt
 \begin{array}[b]{rl}
&\D u\in C([0,\infty);H^N(\mathbb{R}^n)),\ \nabla u\in
L^2([0,\infty),H^{N-1}(\mathbb{R}^n)),\\[2mm]
&\D q\in C([0,\infty);H^N(\mathbb{R}^n))\cap
L^2([0,\infty),H^{N+1}(\mathbb{R}^n)).
\end{array}
$$
Moreover, the solution verifies the following uniform energy
estimate
$$
\|u(t)\|_{H^N(\mathbb{R}^n)}^2+\|q(t)\|_{H^{N+1}(\mathbb{R}^n)}^2+\int_0^t(\|\nabla
u(\tau)\|_{H^{N-1}(\mathbb{R}^n)}^2
+\|q(\tau)\|_{H^{N+1}(\mathbb{R}^n)}^2){\rm d}\tau\leq
CE_0^2.\eqno(1.4)
$$
\bigbreak

Our main results in this paper can be stated as follows:

\vspace{2mm}

\noindent\textbf{Theorem 1.1.} Let $1\leq n\leq4$. Suppose
$\|u_0\|_{H^N(\mathbb{R}^n)}\leq\delta_0\ll1$. If further, $u_0\in
L^1(\mathbb{R}^n)$ (which need not be small), we have
$$
\|D^lu(t)\|_{L^2(\mathbb{R}^n)}\leq
C(1+t)^{-\frac{n}{4}-\frac{l}{2}},\ \ \ \  l=0,1,\cdots,N;
\eqno(1.5)
$$
$$
\|D^lq(t)\|_{L^2(\mathbb{R}^n)}\leq
C(1+t)^{-\frac{n}{4}-\frac{l+1}{2}},\ \ \ \
l=1,2,\cdots,N-1.\eqno(1.6)
$$

\bigbreak \noindent\textbf{Remark 1.1.}  One can rewrite the system
(1.1) as a decouple system of $(u,q)(x,t)$:
$$
\left\{\begin{array}{l}
\D u_t+\sum_{j=1}^nf_j(u)_{x_j}=-u+(I-\Delta)^{-1}u, \\[3mm]
q=-(I-\Delta)^{-1}\nabla u.
 \end{array}
  \right.
 \eqno(1.7)
 $$
To obtain the time-decay rates of $q(x,t)$, from $(1.7)_2$, it
suffices to prove the estimates on $u(x,t)$, i.e. (1.5). When $l=0$,
(1.5) can be obtained directly by using the usual Fourier splitting
method. In addition, for any integer $n\geq 1$, (1.5) also holds
true, see Proposition 3.1 in Section 3. This is consistent with the
result in \cite{Duan1}. When $l\geq1$, to prove (1.5), the estimates
in the Proposition 1.1 and the smallness of $E_0$ should be
employed, see details in the proof of Proposition 3.2.

\bigbreak \noindent\textbf{Theorem 1.2.} Suppose
$\|u_0\|_{H^N(\mathbb{R}^3)}\leq\delta_0\ll1$.  If further, $u_0\in
\dot{H}^{-s}(\mathbb{R}^3)$ (which need not be small), for some
$s\in[0,\frac{3}{2})$, there exists a positive constant $C_0$ such
that
  $$
  \|u(t)\|_{\dot{H}^{-s}(\mathbb{R}^3)}^2\leq C_0,\ \ \ \  \|q(t)\|_{\dot{H}^{-s}(\mathbb{R}^3)}^2\leq C_0, \eqno(1.8)
 $$
and the following decay estimates hold:
$$
\|D^lu(t)\|_{H^{N-l}(\mathbb{R}^3)}\leq C_0(1+t)^{-\frac{l+s}{2}}, \
l=0,\cdots,N, \eqno(1.9)
$$
$$
\|D^lq(t)\|_{H^{N-1-l}(\mathbb{R}^3)}\leq
C_0(1+t)^{-\frac{l+s+1}{2}},\ l=0,\cdots,N-1. \eqno(1.10)
$$

By employing the Hardy-littlewood-Sobolev theorem, for $p\in(1,2]$,
we have $L^p\in\dot{H}^{-s}(\mathbb{R}^3)$ with
$s=3(\frac{1}{p}-\frac{1}{2})\in[0,\frac{3}{2})$. Then, from Theorem
1.2, the following optimal $L^p$-$L^2$ type decay results are
obtained.

\bigbreak \noindent\textbf{Corollary 1.1.($L^p$-$L^2$ time-decay
estimates)} Assume that
$\|u_0\|_{H^N(\mathbb{R}^3)}\leq\delta_0\ll1$. If further,
$u_0(x)\in L^p(\mathbb{R}^3)$ with some $p\in(1,2]$, then the
following decay results hold for any integer $l$ with $0\leq l\leq
N$:
$$
\|D^lu(t)\|_{H^{N-l}(\mathbb{R}^3)}\leq
C_0(1+t)^{-\frac{3}{2}(\frac{1}{p}-\frac{1}{2})-\frac{l}{2}}.\eqno(1.11)
$$
$$
\|D^lq(t)\|_{H^{N+1-l}(\mathbb{R}^3)}\leq
C_0(1+t)^{-\frac{3}{2}(\frac{1}{p}-\frac{1}{2})-\frac{l+1}{2}}.\eqno(1.12)
$$

\bigbreak \noindent\textbf{Remark 1.2.} All the decay results above
are obtained without the smallness of the initial perturbation in
$L^p(\mathbb{R}^3), p\in(1,2]$ or $\dot{H}^{-s}(\mathbb{R}^3)$. The
results generalize those in \cite{Liu1,Wang} for the case of
three-dimensional space. The similar problem in $\mathbb{R}^n$ with
$n=1,2$ will be investigated in future.

\bigbreak \noindent\textbf{Notations.} In this paper, $D^l$ with an
integer $l\geq 0$ stands for the usual any spatial derivatives of
order $l$. For $1\leq p\leq \infty$ and an integer $m\geq 0$, we use
$L^p$ and $W^{m,p}$ denote the usual Lebesgue space
$L^p(\mathbb{R}^n)$ and Sobolev spaces $W^{m,p}(\mathbb{R}^n)$ with
norms $\|\cdot\|_{L^p}$ and $\|\cdot\|_{W^{m,p}}$, respectively, and
set $H^m=W^{m,2}$ with norm $\|\cdot\|_{H^m}$ when $p=2$. In
addition, for $s\in\mathbb{R}$, we define a pseudo-differential
operator $\Lambda^s$ by
$$
\Lambda^sg(x)=\int_{\mathbb{R}^n}|\xi|^s\hat{g}(\xi){\rm
e}^{2\pi\sqrt{-1}x\cdot\xi} {\rm d}\xi,
$$
where $\hat{g}$ denotes the Fourier transform of $g$. We define the
homogeneous Sobolev space $\dot{H}^s$ of all $g$ for which
$\|g\|_{\dot{H}^s}$ is finite, where
$$
\|g\|_{\dot{H}^s}:=\|\Lambda^sg\|_{L^2}=\||\xi|^s\hat{g}\|_{L^2}.
$$
Throughout this paper, we will use a non-positive index $s$. For
convenience, we will change the index to be ``$-s$'' with $s\geq0$.
$C$ denotes a positive generic (generally large) constant that may
vary at different places. The integration domain $\mathbb{R}^3$ will
be always omitted without any ambiguity.

The rest of this paper is arranged as follows. In the next section,
some Sobolev type inequalities and some preliminaries are given for
later use. Section 3 shows the proof of Theorem 1.1 by using Fourier
splitting method. In the last section, we obtain the time-decay
estimates stated in Theorem 1.2.

\section*{ 2.\ Preliminaries}

Firstly, we give some Sobolev inequalities which will be used in the
next two sections.

\bigbreak \noindent\textbf{Lemma 2.1.} (Gagliardo-Nirenberg's
inequality). Let $0\leq m,k\leq l$, then we have
$$
\|D^k g\|_{L^p}\leq C \|D^m
g\|_{L^q}^{1-\theta}\|D^lg\|_{L^r}^\theta
$$
where $k$ satisfies
$$
\frac{1}{p}-\frac{k}{n}=(1-\theta)\left(\frac{1}{q}-\frac{m}{n}\right)+\theta\left(\frac{1}{r}-\frac{l}{n}\right).
$$

\bigbreak \noindent\textbf{Lemma 2.2.} (\cite{Guo}, Lemma A.5) Let
$s\geq0$ and $l\geq0$, then we have
$$
\|D^lg\|_{L^2}\leq
C\|D^{l+1}g\|_{L^2}^{1-\theta}\|g\|_{\dot{H}^{-s}}^\theta,\ {\rm
where}\ \theta=\frac{1}{l+s+1}.
$$

\bigbreak \noindent\textbf{Lemma 2.3.} (\cite{Stein}, Chapter V,
Theorem 1) Let $0<s<n, 1<p<q<\infty,
\frac{1}{q}+\frac{s}{n}=\frac{1}{p}$, then
$$
\|\Lambda^{-s}g\|_{L^q}\leq C\|g\|_{L^p}.
$$
 \bigbreak
Now, when $n=3$, we derive an estimate of  Lyapunov-type which plays
an important role in closing the energy estimates at each $l$-th
level in Section 4.

\vspace{2mm}

 \noindent\textbf{Proposition 2.1.}
 Let $(u,q)(x,t)$ be a solution to the Cauchy problem (1.1)-(1.2) in
 $\mathbb{R}^3$. If the assumptions in Proposition 1.1 hold, we have
$$
\frac{\rm d}{{\rm d}t}\|D^lu(t)\|_{H^1}^2+\|D^{l+1}
u(t)\|_{L^2}^2\leq 0,\ \ \ \ 0\leq l\leq N-1.\eqno(2.1)
$$

\noindent\textbf{Proof.} First, we transform system (1.1) into the
following equivalent decoupled system
$$
\left\{\begin{array}{l}
 \D u_t-\Delta u_t-\Delta u=-(1-\Delta)\sum_{j=1}^nf_j(u)_{x_j}, \ \ \ \  x\in\mathbb{R}^n,\ t>0,\\
  q=-(1-\Delta)^{-1}u, \ \ \ \ \ \ \ \ \ \ \ \ \ \ \  \ \ \ \ \  \qquad\qquad x\in\mathbb{R}^n,\ t>0.
\end{array}\right.\eqno(2.2)
$$

Multiplying $(2.2)_1$ by $u$, then a direct calculation gives
$$\arraycolsep=1.5pt
 \begin{array}[b]{rl}
&\D(u^2+|\nabla u|^2)_t+2|\nabla u|^2-2{\rm
div}\left[u\nabla(u_t+u+\sum_{j=1}^3f_j(u)_{x_j})\right]
+\sum_{j=1}^3\left[2\int_0^uf_j'(\eta)\eta d\eta\right]_{x_j}\\[3mm]
=&\D-\sum_{j=1}^3f_j''(u)u_{x_j}|\nabla u|^2.
\end{array}
$$
Integrating the above equality with respect to $x$ on $\mathbb{R}^3$
and using Lemma 2.1, we get
$$
\frac{\rm d}{{\rm d}t}\|u(t)\|_{H^1}^2+\|D u(t)\|_{L^2}^2\leq
C\|\nabla u(t)\|_{L^\infty}\|D u(t)\|_{L^2}^2.
$$
As a result, we find
$$
\frac{\rm d}{{\rm d}t}\|u(t)\|_{H^1}^2+\|D u(t)\|_{L^2}^2\leq
0.\eqno(2.3)
$$
Since $f_{j}(u)=O(u^2)$ when $u\rightarrow0$, without loss of
generality, let $f_{j}(u)=u^2$. In terms of estimates for the
derivatives of the solution $u$, one can apply $D^l$ on $(2.2)_1$,
and multiply the resulting equality by $D^lu$, then integrating it
with respect to $x$ over $\mathbb{R}^3$, obtaining
$$\arraycolsep=1.5pt
 \begin{array}[b]{rl}
&\D\frac{\rm d}{{\rm d}t}(\|D^lu(t)\|_{L^2}^2+\|\nabla
D^lu(t)\|_{L^2}^2)+2\|\nabla
D^lu(t)\|_{L^2}^2\\[3mm]
=&\D-\sum_{j=1}^3\int D^l(uu_{x_j})D^lu{\rm d}x-\sum_{j=1}^3\int
D^l\nabla(uu_{x_j})D^l\nabla u{\rm d}x =:I_1+I_2.
 \end{array}\eqno(2.4)
$$
For $I_1$, by using H\"{o}lder's inequality and Lemma 2.1, we have
$$\arraycolsep=1.5pt  \begin{array}[b]{rl}
I_1=&\D-\sum_{j=1}^3\int_{\mathbb{R}^3}D^l(uu_{x_j})D^lu{\rm
d}x=-\sum_{j=1}^3\sum_{0\leq k\leq
l}\int_{\mathbb{R}^3}(D^kuD^{l-k}u_{x_j})D^lu{\rm d}x\\[3mm]
\leq &\D C \sum_{0\leq k\leq
l}\|(D^kuD^{l-k+1}u)(t)\|_{L^{\frac{6}{5}}}\|D^lu(t)\|_{L^6}\\[3mm]
\leq &\D C \sum_{0\leq k\leq
l}\|(D^kuD^{l-k+1}u)(t)\|_{L^{\frac{6}{5}}}\|D^{l+1}u(t)\|_{L^2}.
\end{array}
\eqno(2.5)
$$
When $k\leq[\frac{l}{2}]$, by using H\"{o}lder's inequality and
Lemma 2.1, we get
$$\arraycolsep=1.5pt  \begin{array}[b]{rl}
\|(D^kuD^{l-k+1}u)(t)\|_{L^{\frac{6}{5}}}\leq
&\D C\|D^ku(t)\|_{L^3}\|D^{l-k+1}u(t)\|_{L^2}\\[3mm]
\leq &\D
C\|D^mu(t)\|_{L^2}^{1-\frac{k}{l+1}}\|D^{l+1}u(t)\|_{L^2}^{\frac{k}{l+1}}\|u(t)\|_{L^2}^{\frac{k}{l+1}}\|D^{l+1}u(t)\|_{L^2}^{1-\frac{k}{l+1}}\\[3mm]
\leq &\D
C(\|D^mu(t)\|_{L^2}+\|u(t)\|_{L^2})\|D^{l+1}u(t)\|_{L^2}\leq
\frac{1}{12}\|D^{l+1}u(t)\|_{L^2},\end{array} \eqno(2.6)
$$
where we have used the fact $\|u(t)\|_{H^N}\leq CE_0^2\ll1$ in
Proposition 1.1, and $m$ satisfies
$$
\frac{k}{3}-\frac{1}{3}=(\frac{m}{3}-\frac{1}{2})\times(1-\frac{k}{l+1})+(\frac{l+1}{3}-\frac{1}{2})\times\frac{k}{l+1}.\eqno(2.7)
$$
As a result, since $k\leq[\frac{l}{2}]$, we have
$$
m=\frac{l+1}{2(l+1-k)}\in[\frac{1}{2},1).
$$
When $k\geq[\frac{l}{2}]+1$, from H\"{o}lder's inequality and Lemma
2.1 again, we find
$$\arraycolsep=1.5pt  \begin{array}[b]{rl}
\|(D^kuD^{l+1-k}u)(t)\|_{L^{\frac{6}{5}}}\leq
&\D C\|D^ku(t)\|_{L^2}\|D^{l+1-k}u(t)\|_{L^3}\\[3mm]
\leq &\D
C\|u(t)\|_{L^2}^{1-\frac{k}{l+1}}\|D^{l+1}u(t)\|_{L^2}^{\frac{k}{l+1}}\|D^mu(t)\|_{L^2}^{\frac{k}{l+1}}\|D^{l+1}u(t)\|_{L^2}^{1-\frac{k}{l+1}}\\[3mm]
\leq &\D
 C(\|D^mu(t)\|_{L^2}+\|u(t)\|_{L^2})\|D^{l+1}u(t)\|_{L^2}\leq
\frac{1}{4}\|D^{l+1}u(t)\|_{L^2},
\end{array}
\eqno(2.8)
$$
where $m$ is defined by
$$
\frac{l+1-k}{3}-\frac{1}{3}=(\frac{m}{3}-\frac{1}{2})\times\frac{k}{l+1}+(\frac{l+1}{3}-\frac{1}{2})\times(1-\frac{k}{l+1}),
$$
that is, $m=\frac{l+1}{2k}\in(\frac{1}{2},1]$ since
$k\geq\frac{l+1}{2}$. \\
From (2.5), (2.6) and (2.8), we have
$$
I_1=-\sum_{j=1}^3\int D^l(uu_{x_j})D^lu{\rm d}x\leq
\frac{1}{4}\|D^{l+1}u(t)\|_{L^2}.\eqno(2.9)
$$
The estimate of $I_2$ is absolutely the same to (3.10)-(3.12) except
that we replace $l$ by $l+1$. In fact, we have
$$
I_2=-\sum_{j=1}^3\int D^l\nabla(uu_{x_j})D^l\nabla u{\rm d}x\leq
\frac{1}{4}\|D^{l+1}u(t)\|_{L^2}.\eqno(2.10)
$$

Thus, from (2.4), (2.9) and (2.10), we get the estimate of
Lyaponov-type (2.1). Then, we complete the proof of Proposition 2.1.
\qed

\section*{ 3.\ Decay results with initial perturbation in $L^1(\mathbb{R}^n)$}

In this section, we will give optimal decay results by Fourier
splitting method introduced in \cite{Schonbek1,Schonbek2} together
with energy estimates. Theorem 1.1 will be proved by the following
lemmas. First, a straightforward application of Fourier splitting
method yields an optimal $L^2$-norm time-decay rates of solutions as
follows.

\bigbreak \noindent\textbf{Proposition 3.1.} If the initial data
$u_0(x)\in L^1(\mathbb{R}^n)\cap L^2(\mathbb{R}^n)$ with $n\geq1$,
one has
$$
\|u(t)\|_{L^2(\mathbb{R}^n)}^2\leq C
(\|u_0\|_{L^1(\mathbb{R}^n)}+\|u_0\|_{L^2(\mathbb{R}^n)})(1+t)^{-\frac{n}{2}}.
$$
\noindent\textbf{Proof.} First, multiplying $(1.7)_1$ by $u$ and
summing up them and then integrating over $\mathbb{R}^n$, we obtain
$$\arraycolsep=1.5pt  \begin{array}[b]{rl}
\D\frac{1}{2}\frac{\rm d}{{\rm d}t}\int_{\mathbb{R}^n} u^2{\rm d}x
=&\D\int_{\mathbb{R}^n}\sum_{j=1}^nf_j(u)\partial_{x_j}u{\rm
d}x-\int_{\mathbb{R}^n}
u^2{\rm d}x+\int_{\mathbb{R}^n} u(1-\Delta)^{-1}u{\rm d}x\\[3mm]
=&\D-\int_{\mathbb{R}^n} u^2{\rm d}x+\int_{\mathbb{R}^n}
u(1-\Delta)^{-1}u{\rm d}x.
\end{array}
\eqno(3.1)
$$
By Plancherel theorem, we get
$$\arraycolsep=1.5pt  \begin{array}[b]{rl}
&\D\frac{\rm d}{{\rm d}t}\int_{\mathbb{R}^n}|\hat{u}|^2{\rm
d}\xi=-\int_{\mathbb{R}^n}|\hat{u}|^2{\rm d}\xi
+\int_{\mathbb{R}^n}\frac{1}{1+|\xi|^2}|\hat{u}|^2{\rm d}\xi\\[3mm]
= &\D-\int_{\mathbb{R}^n}|\hat{u}|^2{\rm
d}\xi+\int_{|\xi|\leq\sqrt{\frac{n}{t}}}\frac{1}{1+|\xi|^2}|\hat{u}|^2{\rm
d}\xi
+\int_{|\xi|>\sqrt{\frac{n}{t}}}\frac{1}{1+|\xi|^2}|\hat{u}|^2{\rm d}\xi\\[3mm]
\leq &\D-\int_{\mathbb{R}^n}|\hat{u}|^2{\rm
d}\xi+\int_{|\xi|>\sqrt{\frac{n}{t}}}\frac{t}{n+t}|\hat{u}|^2{\rm
d}\xi +\int_{|\xi|
\leq\sqrt{\frac{n}{t}}}\frac{1}{1+|\xi|^2}|\hat{u}|^2{\rm d}\xi\\[3mm]
= &\D-\frac{n}{n+t}\int_{\mathbb{R}^n}|\hat{u}|^2{\rm
d}\xi+\int_{|\xi|\leq\sqrt{\frac{n}{t}}}\frac{1}{1+|\xi|^2}|\hat{u}|^2{\rm
d}\xi
-\frac{t}{n+t}\int_{|\xi|\leq\sqrt{\frac{n}{t}}}|\hat{u}|^2{\rm
d}\xi.
\end{array}
\eqno(3.2)
$$
We rewrite (3.2) as follows.
$$
\frac{\rm d}{{\rm d}t}\int_{\mathbb{R}^n}|\hat{u}|^2{\rm
d}\xi+\frac{n}{n+t}\int_{\mathbb{R}^n}|\hat{u}|^2{\rm
d}\xi\leq\int_{|\xi|\leq\sqrt{\frac{n}{t}}}\left(\frac{1}{1+|\xi|^2}
-\frac{t}{n+t}\right)|\hat{u}|^2{\rm d}\xi,
$$
which is multiplied by $(n+t)^n$ yields
$$\arraycolsep=1.5pt  \begin{array}[b]{rl}
\dis \frac{\rm d}{{\rm
d}t}\left[(n+t)^n\int_{\mathbb{R}^n}|\hat{u}|^2{\rm d}\xi\right]
\leq &\D
(n+t)^n\int_{|\xi|\leq\sqrt{\frac{n}{t}}}\left(\frac{1}{1+|\xi|^2}-\frac{t}{n+t}\right)|\hat{u}|^2{\rm d}\xi\\[3mm]
\leq
&\D(n+t)^n\|\hat{u}(t)\|_{L^\infty(\mathbb{R}_\xi^n)}^2\int_{|\xi|\leq\sqrt{\frac{n}{t}}}\left(1-\frac{t}{n+t}\right){\rm d}\xi\\[3mm]
\leq &\D C\|u(t)\|_{L^1(\mathbb{R}^n)}^2(n+t)^{n-1}(n+t)^{-\frac{n}{2}}\\[3mm]
\leq &\D C\|u_0\|_{L^1(\mathbb{R}^n)}^2(n+t)^{\frac{n}{2}-1},
\end{array}
\eqno(3.3)
$$
where we have used
$\|\hat{u}(t)\|_{L^\infty(\mathbb{R}_\xi^n)}\leq\|u(t)\|_{L^1(\mathbb{R}^n)}$
and the fact $\|u(t)\|_{L^1(\mathbb{R}^n)}\leq
\|u_0\|_{L^1(\mathbb{R}^n)}$ in
\cite{Francesco1,Gao1,Gao2}.\\
Integrating (3.3) with respect to $t$, we have
$$
\|u(t)\|_{L^2(\mathbb{R}^n)}^2\leq
C(\|u_0\|_{L^1(\mathbb{R}^n)}+\|u_0\|_{L^2(\mathbb{R}^n)})(n+t)^{-\frac{n}{2}}.\eqno(3.4)
$$
This proves Proposition 3.1. \qed

\bigbreak

In the proof of Lemma 2.1, we used the essential fact that $\D
\int_{\mathbb{R}^n}\sum_{j=1}^nf_j(u)u_{x_j}{\rm d}x$ in (3.1) is
equal to 0. However, for the derivatives of $u$, the Fourier
splitting method above can not be used directly. The following is
the main reason. Similar to prove (3.1), we have
$$\arraycolsep=1.5pt  \begin{array}[b]{rl}
\D \frac{1}{2}\frac{\rm d}{{\rm d}t}\int_{\mathbb{R}^n}
\!|D^lu|^2{\rm d}x\!
=\!&\D\!\int_{\mathbb{R}^n}\!\sum_{j\!=\!1}^nD^lf_j(u)_{x_j}D^lu{\rm
d}x\!-\!\int_{\mathbb{R}^n} |D^lu|^2{\rm d}x\!+\!\int_{\mathbb{R}^n}
\!D^lu(1\!-\!\Delta)^{\!-\!1}D^lu{\rm d}x.
\end{array}
\eqno(3.5)
$$
But, the term
$\D\int_{\mathbb{R}^n}\sum_{j=1}^nD^lf_j(u)_{x_j}D^lu{\rm d}x$ in
the RHS of (3.5) is not equal to 0. How to control this term? Here
we use the existence results and the smallness of solutions in $H^N$
space stated in Proposition 1.1.

\bigbreak \noindent\textbf{Proposition 3.2.} Let $1\leq n\leq 4$.
Suppose $\|u_0\|_{H^N(\mathbb{R}^n)}\ll1$. If $u_0(x)\in
L^1(\mathbb{R}^n)$ and $D^lu_0(x)\in L^2(\mathbb{R}^n)$, one has
$$
\|D^lu(t)\|_{L^2(\mathbb{R}^n)}^2\leq C
(\|u_0\|_{L^1(\mathbb{R}^n)}+\|D^lu_0\|_{L^2(\mathbb{R}^n)})(1+t)^{-\frac{n}{2}-l},
\ l=1,2,\cdots,N. \eqno(3.6)
$$
\noindent\textbf{Proof.} Applying $D^l\ (1\leq l\leq N)$ on
$(1.7)_1$, and multiplying the resulting equality by $D^lu$, then
integrating it with respect to $x$ over $\mathbb{R}^n$, one has
$$\arraycolsep=1.5pt  \begin{array}[b]{rl}
\D\frac{1}{2}\frac{\rm d}{{\rm d}t}\int_{\mathbb{R}^n} |D^lu|^2{\rm
d}x
\!=\!-\!\int_{\mathbb{R}^n}\sum_{j\!=\!1}^nD^lf_j(u)_{x_j}D^lu{\rm
d}x\!-\!\int_{\mathbb{R}^n}|D^lu|^2{\rm d}x\!+\!\int_{\mathbb{R}^n}
D^lu(1\!-\!\Delta)^{-\!1}D^lu{\rm d}x.
\end{array}
\eqno(3.7)
$$
As mentioned above, to use the Fourier splitting method as
Proposition 3.1, we have to get
$$
\left|\int_{\mathbb{R}^n}\sum_{j=1}^nD^lf_j(u)_{x_j}D^lu{\rm
d}x\right|\leq \delta_0\int_{\mathbb{R}^n} |D^lu|^2{\rm
d}x,\eqno(3.8)
$$
where $0\leq\delta_0\ll1$. In the following, we only prove (3.8) for
$3\leq l\leq N$. For $l=1,2$, the proof of (3.8) is similar. We omit
it here. Without loss of generality, let $f_j(u)=u^2$, one has
$$\arraycolsep=1.5pt  \begin{array}[b]{rl}
&\D\int_{\mathbb{R}^n} D^lf_j(u)_{x_j}D^lu{\rm d}x\\[3mm]
=&\D\int_{\mathbb{R}^n} uD^lu_{x_j}D^lu{\rm d}x +\int_{\mathbb{R}^n}
DuD^{l-1}u_{x_j}D^lu{\rm d}x
+\sum_{2\leq k\leq(l-1)}\int_{\mathbb{R}^n} C_k^lD^k uD^{l+1-k}uD^lu{\rm d}x\\[3mm]
\leq &\D
2\|Du(t)\|_{L^\infty(\mathbb{R}^n)}\|D^lu(t)\|_{L^2(\mathbb{R}^n)}^2
 +C\sum_{2\leq
k\leq(l-1)}\int_{\mathbb{R}^n}
D^k uD^{l+1-k}uD^lu{\rm d}x\\[3mm]
\leq
&\D2\|Du(t)\|_{L^\infty(\mathbb{R}^n)}\|D^lu(t)\|_{L^2(\mathbb{R}^n)}^2
+C\|D^k
u(t)\|_{L^4(\mathbb{R}^n)}\|D^{l-k}u_{x_j}(t)\|_{L^4(\mathbb{R}^n)}\|D^lu(t)\|_{L^2(\mathbb{R}^n)}\\[3mm]
\leq &\D
2\|Du(t)\|_{L^\infty(\mathbb{R}^n)}\|D^lu(t)\|_{L^2(\mathbb{R}^n)}^2\\[3mm]
& \D
+C\|u(t)\|_{L^2(\mathbb{R}^n)}^{1-\theta_1}\|D^lu(t)\|_{L^2(\mathbb{R}^n)}^{\theta_1}\|u(t)\|_{L^2(\mathbb{R}^n)}^{1-\theta_2}\|D^lu(t)\|_{L^2(\mathbb{R}^n)}^{\theta_2}\|D^lu(t)\|_{L^2(\mathbb{R}^n)}\\[3mm]
\leq &\D
2\|Du(t)\|_{L^\infty(\mathbb{R}^n)}\|D^lu(t)\|_{L^2(\mathbb{R}^n)}^2+C\|u(t)\|_{L^2(\mathbb{R}^n)}^{2-(\theta_1+\theta_2)}\|D^lu(t)\|_{L^2(\mathbb{R}^n)}^{\theta_1+\theta_2+1}\\[3mm]
=&\D
2\|Du(t)\|_{L^\infty(\mathbb{R}^n)}\|D^lu(t)\|_{L^2(\mathbb{R}^n)}^2+C\|u(t)\|_{L^2(\mathbb{R}^n)}^{2-(\theta_1+\theta_2)}\|D^lu(t)\|_{L^2(\mathbb{R}^n)}^{\theta_1+\theta_2-1}\|D^lu(t)\|_{L^2(\mathbb{R}^n)}^2,
\end{array}
\eqno(3.9)
$$
where $C_k^l=\left(
                     \begin{array}{c}
                       l \\
                       k \\
                     \end{array}
                   \right)$ and
$$
\theta_1=\frac{k+\frac{n}{4}}{l}\ {\rm and}\
\theta_2=\frac{l-k+\frac{n}{4}+1}{l}.\eqno(3.10)
$$
By noticing $1\leq n\leq 4$ and $2\leq k\leq l-1$, from (3.10), we
know $0<\theta_1\leq1$ and $0<\theta_2\leq1$. Then, it follows from
(3.9) and the fact $\|u(t)\|_{H^N}\ll1$ that
$$
\int_{\mathbb{R}^n} D^lf_j(u)_{x_j}D^lu{\rm d}x\leq
\frac{1}{2}\|D^lu(t)\|_{L^2(\mathbb{R}^n)}^2.\eqno(3.11)
$$
Combining (3.11) and (3.7), we have
$$
\frac{\rm d}{{\rm d}t}\int_{\mathbb{R}^n} (D^lu)^2{\rm d}x \leq
-\int_{\mathbb{R}^n} (D^lu)^2{\rm d}x+\int_{\mathbb{R}^n}
D^lu(1-\Delta)^{-1}D^lu{\rm d}x. \eqno(3.12)
$$
From Plancherel theorem, one has
$$\arraycolsep=1.5pt  \begin{array}[b]{rl}
&\D\frac{\rm d}{{\rm
d}t}\int_{\mathbb{R}^n}\left||\xi|^{2l}\hat{u}\right|^2{\rm d}\xi
\leq-\int_{\mathbb{R}^n}\left||\xi|^{2l}\hat{u}\right|^2{\rm d}\xi
+\int_{\mathbb{R}^n}\frac{1}{1+|\xi|^2}\left||\xi|^{2l}\hat{u}\right|^2{\rm d}\xi\\[3mm]
\leq &\D\!-\!\int_{\mathbb{R}^n}\left||\xi|^{2l}\hat{u}\right|^2{\rm
d}\xi+\int_{|\xi|>\sqrt{\frac{n+2l}{t}}}\frac{t}{n+2l+t}\left||\xi|^{2l}\hat{u}\right|^2{\rm
d}\xi
+\int_{|\xi|\leq\sqrt{\frac{n+2l}{t}}}\frac{1}{1+|\xi|^2}\left||\xi|^{2l}\hat{u}\right|^2{\rm d}\xi\\[3mm]
\leq
&\D\!-\frac{n\!+\!2l}{n\!+\!2l\!+\!t}\int_{\mathbb{R}^n}\left||\xi|^{2l}\hat{u}\right|^2\!{\rm
d}\xi\!+\!\int_{|\xi|\!\leq\!\sqrt{\frac{n\!+\!2l}{t}}}\frac{1}{1\!+\!|\xi|^2}\left||\xi|^{2l}\hat{u}\right|^2\!{\rm
d}\xi
\!-\frac{t}{n\!+\!2l\!+\!t}\int_{|\xi|\!\leq\!\sqrt{\frac{n\!+\!2l}{t}}}\!\left||\xi|^{2l}\hat{u}\right|^2\!{\rm
d}\xi.
\end{array}
\eqno(3.13)
$$
We rewrite (3.13) as
$$\arraycolsep=1.5pt  \begin{array}[b]{rl}
&\D\frac{\rm d}{{\rm
d}t}\int_{\mathbb{R}^n}\left||\xi|^{2l}\hat{u}\right|^2{\rm
d}\xi+\frac{n+2l}{n+2l+t}\int_{\mathbb{R}^n}\left||\xi|^{2l}\hat{u}\right|^2{\rm
d}\xi\\[3mm]
\leq &\D\int_{|\xi|
\leq\sqrt{\frac{n+2l}{t}}}\left(\frac{1}{1+|\xi|^2}
-\frac{t}{n+2l+t}\right)\left||\xi|^{2l}\hat{u}\right|^2{\rm d}\xi.
\end{array}
\eqno(3.14)
$$
Consequently, multiplying (3.14) by $(n+2l+t)^{n+2l}$, one has
$$\arraycolsep=1.5pt  \begin{array}[b]{rl}
&\D\frac{\rm d}{{\rm d}t}\left[(n+2l+t)^{n+2l}\int_{\mathbb{R}^n}\left||\xi|^{2l}\hat{u}\right|^2{\rm d}\xi\right]\\[3mm]
\leq
&\D(n+2l+t)^{n+2l}\int_{|\xi|\leq\sqrt{\frac{n+2l}{t}}}\left(\frac{1}{1+|\xi|^2}-\frac{t}{n+2l+t}\right)\left||\xi|^{2l}\hat{u}\right|^2{\rm d}\xi\\[3mm]
\leq
&\D(n+2l+t)^{n+2l}\|\hat{u}(t)\|_{L^\infty(\mathbb{R}_\xi^n)}^2\int_{|\xi|\leq\sqrt{\frac{n+2l}{t}}}\left(1-\frac{t}{n+2l+t}\right)|\xi|^{2l}{\rm d}\xi\\[3mm]
\leq &\D C\|u(t)\|_{L^1(\mathbb{R}^n)}^2(n+2l+t)^{n+2l-1}(n+2l+t)^{-\frac{n+2l}{2}}\\[3mm]
\leq &\D C\|u_0\|_{L^1(\mathbb{R}^n)}^2(n+2l+t)^{\frac{n+2l}{2}-1}.
\end{array}
\eqno(3.15)
$$
Integrating (3.15) with respect to $t$ over $(0,t)$, we find
$$
\|D^lu(t)\|_{L^2(\mathbb{R}^n)}^2=\||\xi|^l\hat{u}(t)\|_{L^2(\mathbb{R}_\xi^n)}^2\leq
C(\|u_0\|_{L^1(\mathbb{R}^n)}+\|D^lu_0\|_{L^2(\mathbb{R}^n)})(n+2l+t)^{-\frac{n+2l}{2}}.\eqno(3.16)
$$
Then, we get (3.6). \qed

\section*{ 4.\ Decay results with initial perturbation in $\dot{H}^{-s}(\mathbb{R}^3)$}

This section devotes to the optimal $L^p$-$L^2(\mathbb{R}^3)$ decay
rates of solutions to (1.1)-(1.2) when the initial data is in the
negative Sobolev space $\dot{H}^{-s}(\mathbb{R}^3)$ with
$s\in[0,\frac{3}{2})$.

The following lemma plays a key role in the proof of Theorem 1.2. It
shows an energy estimate of the solutions in the negative Sobolev
space $\dot{H}^{-s}(\mathbb{R}^3)$. Namely, we have

 \noindent\textbf{Lemma 4.1.} If
$\mathcal{E}_0:=\|u_0\|_{H^N}\ll1$, for $s\in(0,\frac{1}{2}]$, we
have
$$\arraycolsep=1.5pt  \begin{array}[b]{rl}
 &\D \frac{\rm d}{{\rm d}t}\int\left(|\Lambda^{-s}u|^2+|\Lambda^{-s}\nabla
 u|^2\right){\rm d}x+\int|\nabla\Lambda^{-s}u|^2{\rm d}x\\[3mm]
 \leq &\D C(\|Du(t)\|_{H^1}^2+\|D^2u(t)\|_{H^1}^2)(\|\Lambda^{-s}u(t)\|_{L^2}+\|\Lambda^{-s}\nabla
 u(t)\|_{L^2});
 \end{array}
 \eqno(4.1)
$$
and for $s\in(\frac{1}{2},\frac{3}{2})$, we have
$$
\arraycolsep=1.5pt  \begin{array}[b]{rl} &\D \frac{\rm d}{{\rm
d}t}\int\left(|\Lambda^{-s}u|^2+|\Lambda^{-s}\nabla
u|^2\right){\rm d}x+\int|\nabla\Lambda^{-s}u|^2{\rm d}x\\[3mm]
\leq &\D
C(\|u(t)\|_{L^2}^{s-\frac{1}{2}}\|Du(t)\|_{L^2}^{\frac{5}{2}-s}+\|Du(t)\|_{L^2}^{s-\frac{1}{2}}\|D^2u(t)\|_{L^2}^{\frac{5}{2}-s})
(\|\Lambda^{-s}u(t)\|_{L^2}+\|\Lambda^{-s}\nabla u(t)\|_{L^2}).
\end{array}
\eqno(4.2)
$$
\noindent\textbf{Proof.} Applying $\Lambda^{-s}$ to $(2.2)_1$ and
multiplying the resulting identity by $\Lambda^{-s}u$, and
integrating over $\mathbb{R}^3$ by parts, we get
$$\arraycolsep=1.5pt  \begin{array}[b]{rl}
&\D\frac{1}{2}\frac{\rm d}{{\rm
d}t}\int\left(|\Lambda^{-s}u|^2+|\Lambda^{-s}\nabla
u|^2\right){\rm d}x+\int|\nabla\Lambda^{-s}u|^2{\rm d}x\\[3mm]
=&\D-\sum_{j=1}^3\left\{\int\Lambda^{-s}u\Lambda^{-s}f_j(u)_{x_j}{\rm d}x-\int\Lambda^{-s}u\Delta\Lambda^{-s}f_j(u)_{x_j}{\rm d}x\right\}\\[3mm]
:=&\D J_1+J_2.
\end{array}
\eqno(4.3)
$$
For $J_1$, using H\"{o}lder inequality, Lemma 2.1, Lemma 2.3 and
Young's inequality, we have
$$\arraycolsep=1.5pt  \begin{array}[b]{rl}
J_1\leq &\D C\|\Lambda^{-s}u(t)\|_{L^2}\sum_{j=1}^3\|\Lambda^{-s}f_j(u)_{x_j}(t)\|_{L^2}\\[3mm]
\leq &\D
C\|\Lambda^{-s}u(t)\|_{L^2}\sum_{j=1}^3\|f_j(u)_{x_j}(t)\|_{L^{\frac{1}{\frac{1}{2}+\frac{s}{3}}}}
\leq  C\|\Lambda^{-s}u(t)\|_{L^2}\|\nabla
u(t)\|_{L^2}\|u(t)\|_{L^{\frac{3}{s}}}\\[3mm]
\leq &\D C\|\Lambda^{-s}u(t)\|_{L^2}\|\nabla u(t)\|_{L^2}\|\nabla
u(t)\|_{L^2}^{\frac{1}{2}+s}\|D^2u(t)\|_{L^2}^{\frac{1}{2}-s}\\[3mm]
\leq &\D C(\|Du(t)\|_{L^2}^2+\|D^2u(t)\|_{L^2}^2)\|\Lambda^{-s}
u(t)\|_{L^2}.
\end{array}
\eqno(4.4)
$$
Here we have used the facts $\frac{1}{2}+\frac{s}{3}<1$ and
$\frac{3}{s}\geq6$. \\
Similarly, it holds that
$$
J_2\leq
C(\|D^2u(t)\|_{L^2}^2+\|D^3u(t)\|_{L^2}^2)\|\Lambda^{-s}\nabla
u(t)\|_{L^2}.\eqno(4.5)
$$
Combining (4.4) and (4.5), we get (4.1).\\
Next, we want to prove (4.2). A direct calculation as (4.4) gives
$$
\arraycolsep=1.5pt  \begin{array}[b]{rl} J_1
\leq &\D C\|\Lambda^{-s}u(t)\|_{L^2}\sum_{j=1}^3\|\Lambda^{-s}f_j(u)_{x_j}(t)\|_{L^2}\\[3mm]
\leq &\D C\|\Lambda^{-s}u(t)\|_{L^2}\sum_{j=1}^3\|f_j(u)_{x_j}(t)\|_{L^{\frac{1}{\frac{1}{2}+\frac{s}{3}}}}\\[3mm]
\leq &\D C\|\Lambda^{-s}u(t)\|_{L^2}\|D
u(t)\|_{L^2}\|u(t)\|_{L^{\frac{3}{s}}} \leq
C\|\Lambda^{-s}u(t)\|_{L^2}\|Du(t)\|_{L^2}
\|u(t)\|_{L^2}^{s-\frac{1}{2}}\|Du(t)\|_{L^2}^{\frac{3}{2}-s}\\[3mm]
\leq &\D
C\|u(t)\|_{L^2}^{s-\frac{1}{2}}\|Du(t)\|_{L^2}^{\frac{5}{2}-s}\|\Lambda^{-s}
u(t)\|_{L^2}.
\end{array}
\eqno(4.6)
$$
In the same way, one has
$$
I_2\leq
C\|Du(t)\|_{L^2}^{s-\frac{1}{2}}\|D^2u(t)\|_{L^2}^{\frac{5}{2}-s}\|\Lambda^{-s}
u(t)\|_{L^2}.
$$
Thus we complete the proof of Lemma 4.1.  \qed

With Proposition 1.1, Proposition 2.1 and Lemma 4.1 in hand, we are
now ready to prove Theorem 1.2.

\noindent\emph{Proof of Theorem 1.2 for the case of
$s\in(0,\frac{1}{2}]$}. Define
$\mathcal{E}_{-s}(t)=\|\Lambda^{-s}u(t)\|_{L^2}^2+\|\Lambda^{-s}\nabla
u(t)\|_{L^2}^2$. Integrating (4.1) with respect to $t$, we find for
$s\in(0,\frac{1}{2}]$,
$$
\mathcal{E}_{-s}(t)\leq
\mathcal{E}_{-s}(0)+C\int_0^t\|Du(\tau)\|_{H^2}^2\sqrt{\mathcal{E}_{-s}(\tau)}{\rm
d}\tau.\eqno(4.7)
$$
From (2.1), we have the estimates of integrability of
$\|Du\|_{H^1}^2$ with respect to $t$. As a result, we have
$$
\mathcal{E}_{-s}(t)\leq
\mathcal{E}_{-s}(0)+C\sup\limits_{0\leq\tau\leq
t}\sqrt{\mathcal{E}_{-s}(\tau)}\leq C(1+\sup\limits_{0\leq\tau\leq
t}\sqrt{\mathcal{E}_{-s}(\tau)}).
$$
This yields $\mathcal{E}_{-s}(t)\leq C_1$ with a positive constant
$C_1$, that is
$$
\|\Lambda^{-s}u(t)\|_{L^2}^2+\|\Lambda^{-s}q(t)\|_{L^2}^2\leq\|\Lambda^{-s}u(t)\|_{L^2}^2+\|\Lambda^{-s}\nabla
u(t)\|_{L^2}^2\leq C_1.\eqno(4.8)
$$
This proves (1.8) for $s\in(0,\frac{1}{2}]$. \\
Next, we recall the estimate of Lyaponov-type (2.1) as
$$
\frac{\rm d}{{\rm d}t}\|D^lu(t)\|_{H^1}^2+\|D^{l+1}
u(t)\|_{L^2}^2\leq 0,\ \ \ \  0\leq l\leq N-1.
$$
We may use Lemma 2.2 to have
$$
\|D^{l+1}u(t)\|_{L^2}\geq
C\|\Lambda^{-s}u(t)\|_{L^2}^{-\frac{1}{l+s}}\|D^lu(t)\|_{L^2}^{1+\frac{1}{l+s}}.
$$
From the above inequality and (4.8) we get for  each $l$ with $0\leq
l\leq N-1$,
$$\arraycolsep=1.5pt
 \begin{array}[b]{rl}
\|D^{l+1}u(t)\|_{L^2}^2\geq &\D\frac{1}{2}\|D^{l+1}u(t)\|_{L^2}^2+\frac{1}{2}CC_1^{-\frac{1}{l+s}}\|D^lu(t)\|_{L^2}^{2(1+\frac{1}{l+s})}\\[3mm]
\geq &\D
C_2(\|D^{l+1}u(t)\|_{L^2}^2+\|D^lu(t)\|_{L^2}^2)^{1+\frac{1}{l+s}},
\end{array}
\eqno(4.9)
$$
where $C_2>0$ is a constant. Thus we deduce the following time
differential inequality
$$
\frac{\rm d}{{\rm
d}t}\|D^lu(t)\|_{H^{N-l}}^2+C_2(\|D^lu(t)\|_{H^{N-l}}^2)^{1+\frac{1}{l+s}}\leq0,\
{\rm for}\ l=0,1,\cdots,N.
$$
Integrating this inequality, one gets for some constant $C_3>0$
$$
\|D^lu(t)\|_{H^{N-l}}^2\leq C_3(1+t)^{-(l+s)},\ {\rm for}\
l=0,1,\cdots,N.\eqno(4.10)
$$
$(2.2)_2$ and (4.10) yield
$$
 \|D^lq(t)\|_{H^{N-1-l}}\leq C_4(1+t)^{-\frac{l+s+1}{2}},\ {\rm for}\
l=0,\cdots,N-1. \eqno(4.11)
$$

\noindent\emph{Proof of Theorem 1.2 for the case of
$s\in(\frac{1}{2},\frac{3}{2})$}. First we can give from what we
have proved for (1.8)-(1.9) with $s=\frac{1}{2}$ above that the
following decay results holds:
$$
\|D^lu(t)\|_{H^{N-l}}^2\leq C_3(1+t)^{-l-\frac{1}{2}},\ {\rm for}\
l=0,1,\cdots,N.\eqno(4.12)
$$
As a result, and using (4.2) that for
$s\in(\frac{1}{2},\frac{3}{2})$,
$$\arraycolsep=1.5pt
 \begin{array}[b]{rl}
\mathcal{E}_{-s}(t)\leq
&\D\mathcal{E}_{-s}(0)+C\int_0^t\|u(\tau)\|_{L^2}^{s-\frac{1}{2}}\|Du(\tau)\|_{L^2}^{\frac{5}{2}-s}\sqrt{\mathcal{E}_{-s}(\tau)}{\rm d}\tau\\[3mm]
\leq &\D C_1+CC_3\int_0^t(1+\tau)^{-\frac{7}{4}-\frac{s}{2}}{\rm
d}\tau\cdot\sup\limits_{0\leq\tau\leq
t}\sqrt{\mathcal{E}_{-s}(\tau)}\\[3mm]
\leq &\D C_5(1+\sup\limits_{0\leq\tau\leq
t}\sqrt{\mathcal{E}_{-s}(\tau)}),
 \end{array}
 \eqno(4.13)
$$
which yields (1.8) with $s\in(\frac{1}{2},\frac{3}{2})$. At last,
the proof of (1.9) with $s\in(\frac{1}{2},\frac{3}{2})$ can be
treated as the case of $s\in(0,\frac{1}{2}]$ above.

 Thus, by
taking $C_0=\max\limits_{1\leq i\leq 5}\{C_i\}$, we complete the
proof of Theorem 1.2.
\\
\\
{\bf Acknowledgement:} \ \ The research of Zhigang Wu was supported
by NSFC (No. 11101112) and in part by NSFC (No. 11071162). The
research of Wenjun Wang was supported by the Tian Yuan Fund of
Mathematics in China (No.11126096), NSFC (No. 11201300), Shanghai
university young teacher training program (No. slg11032) and in part
by NSFC (No.11171220, No.11171212).

\bibliographystyle{plain}

\end{document}